\documentclass[10pt]{article}

\usepackage{graphicx}

\date{}
\usepackage[latin1]{inputenc}
\usepackage[english]{babel}
\usepackage[T1]{fontenc}
\usepackage{amssymb}
\usepackage{amsthm}    
\usepackage{stmaryrd}
\usepackage{amsmath,amsfonts}
\usepackage{graphicx}
\usepackage{epstopdf}
\DeclareGraphicsExtensions{.pdf,.eps,.png,.jpg,.mps}
\usepackage{indentfirst}
\usepackage{dsfont}
\usepackage{color}
\usepackage{enumitem}
\usepackage{float}
\usepackage{pict2e}
\usepackage{soul}

\newtheorem{theorem}{Theorem}[section]
\newtheorem{lemma}[theorem]{Lemma}

\newcommand{\sgn}{\operatorname{sgn}}

\title{The Mean/Max Statistic in Extreme Value Analysis}
\author{Paul Rochet and Isabel Serra}
\date{}

\begin{document}
\maketitle


\begin{abstract} Most extreme events in real life can be faithfully modeled as random realizations from a Generalized Pareto distribution, which depends on two parameters: the scale and the shape. In many actual situations, one is mostly concerned with the shape parameter, also called tail index, as it contains the main information on the likelihood of extreme events. In this paper, we show that the mean/max statistic, that is the empirical mean divided by the maximal value of the sample, constitutes an ideal normalization to study the tail index independently of the scale. This statistic appears naturally when trying to distinguish between uniform and exponential distributions, the two transitional phases of the Generalized Pareto model. We propose a simple methodology based on the mean/max statistic to detect, classify and infer on the tail of the distribution of a sample. Applications to seismic events and detection of saturation in experimental measurements are presented.
\end{abstract}



\section{Introduction}
Two fundamental results marked the starting point of Extreme Value Theory (EVT): Fisher-Tippet-Gnedenko and Pickands-Balkemax-de Haan Theorems. These theorems enabled to characterize the asymptotic behavior of extreme values by a real number $k$, called tail index (we refer to \cite{markovich}, \cite{beirlant2006statistics}, \cite{coles2001introduction,embrechts1997modelling,finkenstadt2003extreme} and \cite{mcneil2010quantitative} for a survey). Ever since Gnedenko's original paper, estimation of the tail index has become a major concern in the literature, for which even the most used techniques can produce unsatisfactory results in many situations. For instance, several authors have highlighted the estimation problems that arise from the second approach, see \cite{castillo,hosking1987parameter,zhang2009new}.\\

Pickands-Balkemax-de Haan Theorem has widened the use of the Generalized Pareto Distribution (GPD) in extreme values theory as a model for tails. According to the behavior of the probability density functions of the GPD, we can distinguish three submodels corresponding to $k<0$, $0<k<1$ and $k > 1$ respectively, separated by the exponential distribution ($k=0$) and the uniform distribution ($k=1$). Since no method is able to provide a satisfying estimation of the tail index regardless of the true model (see \cite{del2015likelihood}), the parameter inference should be made subsequently to the choice of submodel, even more so for small samples. The first recommendation is to specify the submodel of the GPD a priori, in order to estimate the parameters.
Therefore, we propose to use a test procedure to distinguish between the two transitional phases of the tail index in GPD models: the exponential and uniform distributions. Due to the monotonic behavior of the test statistic with respect to the tail index of the GPD, the test can be used for classification purposes.\\

Tests on separate families of hypotheses have been originally studied in \cite{cox1961tests} who proposed a generalization of the Neyman-Pearson principle, based on a likelihood ratio criterion. The so-called maximum-likelihood ratio test is not necessarily optimal, although it provides a practical methodology to test a wide range of hypotheses. Examples studied in the literature include invariant tests of exponential versus normal or uniform distributions in \cite{uthoff1970optimum}, normal versus Cauchy distribution in \cite{franck1981most} and so on.\\

In this paper, we show that Cox's maximum likelihood ratio test for uniform versus exponential distributions is the most powerful among scale-free tests. It relies on the ratio $\tau_n$ between
the empirical mean and sample maximum, whose behavior is highly dictated by the tail of the distribution, thus making it a main feature of the GPD framework. The mean/max statistic $\tau_n$ then provides a simple yet effective way to detect the tail behavior in a Generalized Pareto model. In fact, the tail classification method based solely on the value of $\tau_n$ produces extremely conclusive results that compare favorably with the more computationally expensive inference methods such as Zhang and Stephens \cite{zhang2009new}, Song and Song \cite{song2012quantile} or maximum likelihood.\\

The problem of bad classification in GPD estimation was recently detected \cite{del2015likelihood} but no satisfactory solution has been proposed so far. Summarizing, the quality of inference methods for the tail index $k$ in the GPD highly depends on the underlying submode, whether it is Model A ($k >1$), Model B ($0\leq k \leq 1$) or Model C ($k <0$). Because standard methods tend to specialize to certain regions of the values of the parameter, it is crucial to separate the GPD model into these three submodels reflecting the different behavior of the tails.  \\

The test on separate families of hypotheses to distinguish between uniform and exponential distributions is described in Section \ref{sec:test}. In Section \ref{sec:classif}, we show how to extend the test procedure for tail classification and detection in Generalized Pareto models or more general distributions. In Section \ref{sec:appli}, we propose a simple general recipe to infer on the tail of the distribution of a sample and an application on global seismic activity data is presented. Theoretical results regarding the distribution of the mean/max statistic are discussed in the Appendix.

%
%

\section{Test to distinguish between uniform and exponential tail distributions}\label{sec:test}

Let $X=(X_1,\dots ,X_n)$ be a sample of independent identically distributed variables drawn from a distribution on $[0,+ \infty)$ assumed to be either uniform on some interval $[0,\theta]$ with $\theta >0$ or exponential with parameter $\lambda >0$. As it is usually the case when dealing with separate families of hypotheses, a uniformly most powerful test does not exist because the likelihood ratio statistic depends on the true value of the parameter. The maximum likelihood ratio statistic was proposed in \cite{cox1961tests} as a way to generalize Neyman-Pearson's principle to families of distributions, although the optimality of the test is no longer backed up by the theory. Nevertheless, we show that in this situation, Cox's maximum likelihood ratio test is the most powerful among a large class of tests whose distribution is invariant within each hypothesis. For such tests, both the level and power are determined by the common distributions of the statistic under the null hypothesis and the alternative. The unicity of these distributions within each hypothesis enables to deduce the most powerful test from a direct application of Neyman-Pearson's lemma. \\

\noindent We are interested in testing the null hypothesis\vspace{0.2cm}

$ \mathcal H_0: "\text{The } X_i \text{'s are uniformly distributed on an interval} \ [0,\theta] \ \text{with $\theta>0$}"$ \vspace{0.2cm}

\noindent against the alternative \vspace{0.2cm}

$  \mathcal H_1: "\text{The } X_i \text{'s are exponentially distributed.}"  $\vspace{0.2cm}

\noindent In this particular situation, one can exploit the fact that both the uniform and exponential models are closed by positive scaling: if the sample $X$ belongs to the uniform (resp. exponential) model, then so does $t X$ for all $t>0$. We consider the class of scale-free tests, that is binary valued functions $\Phi = \Phi(X)$ of the sample $X$ satisfying
$$ \Phi(X) = \Phi( t X) \ , \ \forall t >0.   $$
If we let $\Phi(X) = \mathds 1 \{ X \in \mathcal R \}$ for $\mathcal R$ a Borel set of $\mathbb R^n$, saying that $\Phi$ is scale-free simply means that $\mathcal R$ is a cone. Equivalently, a statistic is scale-free if it can be expressed as a function of the normalized sample $X/ \Vert X \Vert$, for any norm $\Vert . \Vert$ on $\mathbb R^n$. \\

Because the uniform and exponential models are stable by positive scaling, the distribution of a scale-free statistic does not depend on the parameter, whether we are under the null hypothesis (uniform) or the alternative (exponential). Therefore, we know there exists a most powerful scale-free test for these hypotheses, which we can derive from the likelihood ratio statistic of the normalized sample. As we show below, the resulting test is function of the statistic
$$ \tau_n  := \frac{\overline{X}_n}{X_{(n)}}  $$
where $X_{(n)} = \max \{ X_1,\dots ,X_n \}$ and $\overline{X}_n = \frac 1 n \sum_{i=1}^n X_i$, whose properties are often discussed in relations with the tail distribution, see \cite{villasenor2009bootstrap}.

\begin{theorem}\label{th1} The test $ \Phi = \mathds 1 \{ \tau_n \leq c_\alpha \}$ with $c_\alpha$ such that $\mathbb P (\tau_n \leq c_\alpha \vert \mathcal H_0 ) = \alpha$ is the most powerful scale-free test of level $\alpha \in (0,1)$ to test $\mathcal H_0$ against $\mathcal H_1$.
\end{theorem}

\noindent \textit{Proof.} Let $\Phi = \Phi(X)$ be a scale-free test of level $\alpha \in (0,1)$, we can write almost surely (for $X_n \neq 0$)
$$ \Phi (X_1,\dots ,X_n) = \Phi \Big( \frac{ X_1 }{X_n},\dots , \frac{ X_{n-1} }{X_n},1 \Big) := \Psi \Big( \frac{ X_1 }{X_n},\dots , \frac{ X_{n-1} }{X_n}\Big),$$
for $\Psi: \mathbb R^{n-1} \to \{ 0,1\}$. Let $Y = (Y_1,\dots ,Y_{n-1}) = (X_1 /X_n, \dots  , X_{n-1}/X_n)$. Since the distribution of $Y$ is constant over $\mathcal H_0$ and over $\mathcal H_1$, Neyman-Pearson's Lemma tells us that the likelihood ratio test has maximal power among tests of significance level $\alpha$. Both likelihood functions $\mathcal L_0$ and $\mathcal L_1$ of $Y$ under $\mathcal H_0$ and $\mathcal H_1$ respectively can be calculated explicitly, yielding
$$ \mathcal L_0(Y) = \frac{1}{n (\max \{ 1,Y_{(n-1)} \})^n} \ \ \text{ and } \ \ \mathcal L_1(Y) =  \frac{(n-1)! }{(1+ \sum_{i=1}^{n-1} Y_i)^n}.  $$
Recall that $Y_i = X_i/X_n$, Neyman-Pearson's likelihood ratio is thus given by
$$  \frac{1}{n!} \frac {(1+ \sum_{i=1}^{n-1} Y_i)^n}{ (\max \{ 1,Y_{(n-1)} \})^n}  = \frac{(n \tau_n)^n}{n!} . $$
Since it is an increasing function of $\tau_n$, the likelihood ratio test writes as $  \mathds 1 \{ \tau_n \leq c_\alpha \} $ for some suitable $c_\alpha$.
\qed \\

In this framework, the most powerful scale-free test actually recovers Cox's maximum likelihood ratio test, introduced in \cite{cox1961tests} in a more general context. Indeed, considering the likelihood functions $L_0(\theta, X) = \theta^{-n} \mathds 1 \{ X_i \leq \theta,  i=1,\dots ,n \}$ and $L_1(\lambda, X) = \lambda^{n} \exp(- \lambda \sum_{i=1}^n X_i)$, the maximum likelihood ratio statistic is given by
$$ \frac{\sup_{\lambda >0} L_1(\lambda, X)} {\sup_{\theta >0} L_0(\theta, X)} = \left(\frac{X_{(n)}}{e\ \overline{X}_n}  \right)^n = \left(\frac{1}{e \ \tau_n}  \right)^n.  $$
Thus, Cox's maximum likelihood ratio test rejects the null hypothesis $\mathcal H_0$ for sufficiently small values of $\tau_n$. The threshold $c_\alpha$ corresponding to a test of significance level $\alpha$ can be computed easily by Monte-Carlo.  Nevertheless, to get an analytical expression of the threshold and the corresponding power of the test, one needs to derive the true distribution of the statistic $\tau_n$ under the null hypothesis and under the alternative. This issue is discussed in the Appendix.

\section{Tail classification in Generalized Pareto distribution}\label{sec:classif}

If there is no particular reason to favor the uniform model over the exponential one, the significance level can be chosen so that the probabilities of error under $\mathcal H_0$ and $\mathcal H_1$ are equal, thus inducing a minimal probability of error in the worst case scenario. In this purpose, one must choose the threshold $c_n$ as the unique solution of
$$ \mathbb P(\tau_n \geq c_n \vert \mathcal H_0) = \mathbb P( \tau_n \leq c_n \vert \mathcal H_1) := A(n).   $$
Thus, when using the threshold $c_n$, the probability $A(n)$ of selecting the correct model between uniform and exponential no longer depends on the true distribution. The values of the threshold $c_n$ and percentage of accurate selections $A(n)$ are computed in Table \ref{tab:accuracy} for different sample sizes $n$. It appears that the test quickly reaches a near perfect accuracy as the sample size increases. For a sample of size $n=5$, the method selects the correct model more than $71 \%$ of the time, while a sample of size $n=20$ achieves a precision above $95 \%$.\\

\begin{table}[h]
\begin{center}
\begin{tabular}{|c||c|c|c|c|c|}
\hline $n$ & $5$ & $10$ & $20$ & $50$ & $100$ \\ \hline \hline
        $A(n)$ & $0.7178$ & $0.8437$ & $0.9514$ & $0.9984$ & $1$  \\ \hline $c_n$ & $0.531$ & $0.460$ & $0.420$ & $0.392$ & $0.380$  \\ \hline
\end{tabular}
\caption{\footnotesize{Optimal threshold $c_n$ and percentage of accuracy $A(n)$ for sample sizes $n=5,10,20,50,100$, computed by Monte-Carlo simulations with $5.10^4$ replications.}}
\label{tab:accuracy}
\end{center}
\end{table}

\begin{figure}[h!]
    \centering
        \includegraphics[width=0.99\textwidth]{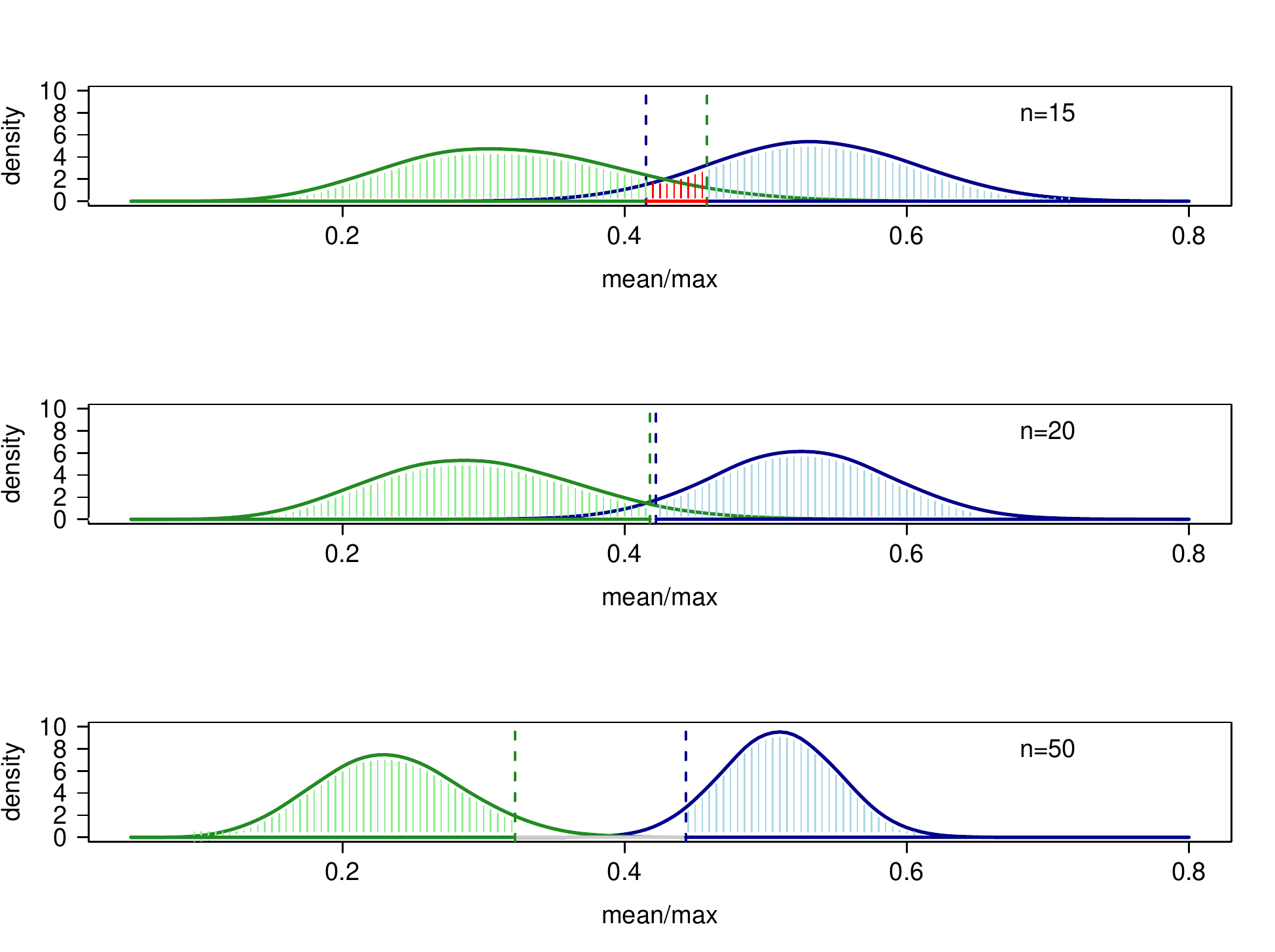}
    \label{fig:ex_test_sizes}
    \caption{\footnotesize{Monte-Carlo estimated densities of the mean/max statistic $\tau_n$ under $\mathcal H_0$ (blue) and $\mathcal H_1$ (green) for sample sizes $15,20$ and $50$ from top to bottom. The  $5\%$ quantile under $\mathcal H_0$ and $95\%$ quantile under $\mathcal H_1$ are marked in vertical dotted lines, pointing out that a sample size of $n=20$ is needed for a $95 \%$ precision. }}
\end{figure}

Of course, the level of accuracy is only exact if the model is well-specified, which is rarely the case in practice. Nevertheless, the test procedure may provide substantial information on the tail of the distribution. To study the behavior of extreme events, one focuses generally on the tail observations, i.e. the data over a certain threshold $\mu$, so as to obtain independent realizations $X_i$ conditionally to $X_i > \mu$. Under regularity conditions, the tail observations (translated to the left by a factor $\mu$ so as to have their support starting at zero) converge in distribution as $\mu \to \infty$ towards independent realizations of a Generalized Pareto Distribution (GPD), with density
$$ f_{k,\sigma}(x) = \frac 1 \sigma \Big(1-k \frac x \sigma \Big)^{\frac{1-k}{k}} \ , \ \left\{ \begin{array}{ll} x \geq 0 & \text{if } k \leq 0 \\ 0 \leq x \leq \frac \sigma k & \text{if } k >0. \end{array} \right. $$
The scale $\sigma >0$ and the shape $k \in \mathbb R$ (also called the tail index) are the two parameters used to describe the tail. The Generalized Pareto density can be extended by continuity to the values $k=0$ and $k=1$ recovering respectively the exponential and uniform distributions.  Generally, the relevant information on the tail comes from the scale parameter $k$ and less importance is given to the scale $\sigma$.\\

When working with tail observations, it is thus customary to assume that the data are drawn from a GPD. In this model, the behavior of the mean/max statistic $\tau_n$ is highly dictated by the tail index $k$ while its distribution does not depend on the scale parameter $\sigma$ ($\tau_n$ is scale-free). The mean/max statistic is thus perfectly suited to investigate the tail index of the GPD without having to concern about the scale. As a matter of fact, a crucial matter in Generalized Pareto models is to determine to what submodel the tail belongs, i.e. if the tail index $k$ is greater than $1$ (model A), between $0$ and $1$ (model B) or negative (model C). Because the frontiers between the different models are achieved for the exponential ($k=0$) and uniform ($k=1$) cases, the test procedure to distinguish between uniform and exponential distributions can be naturally extended for classification purposes. \\

The asymptotic distribution of $\tau_n$ (see Theorem \ref{gpd_asymp} in the Appendix) reveals that $\tau_n$ vanishes at a rate of $n^{-k}$ when $k \in (0,1)$, converges to $k/(k+1)$ when $k>0$ and is of the order $1/\log(n)$ in-between, for $k=0$. In particular, the median of $\tau_n$ can be deduced from Theorem \ref{gpd_asymp},
\begin{equation}\label{med_asymp} \text{med}(\tau_n) = \left\{ \begin{array}{ll}
- \dfrac{k}{k+1} \Big( \dfrac{\log(2)}{n} \Big)^{-k} + o (n^k) & \text{ for } -1 < k <0 \\
& \\
\dfrac{1}{\log(n)} + \dfrac{\log \log(2)}{\log ^2(n)} + o \Big( \dfrac 1 {\log ^2(n)} \Big) &  \text{ for } k =0 \\
& \\
\dfrac{k}{k+1} + o \Big( \dfrac 1 {\sqrt n} \Big) & \text{ for } k >0. \end{array} \right. \end{equation}

The monotonic behavior of $\tau_n$ with respect to the tail index $k$ makes it a usefool tool to learn to which submodel the distribution of the data belongs. The idea is simple: the practitioner chooses two reals numbers $0 < a_n < b_n < 1$ and concludes to the model A if $\tau_n \geq b_n$, the model B if $a_n < \tau_n < b_n$ and the model C if $\tau_n \leq a_n$. As suitable values of $a_n$ and $b_n$, we use the theoretical median of $\tau_n$ under the exponential and uniform distributions respectively. By taking these values, the probability of selecting each model in the transitional cases $k=0$ and $k=1$ does not exceed $1/2$ so that none of the model A, B or C is favored. The actual value of $\text{med}(\tau_n)$ in the uniform case follows directly from Lemma \ref{ih}. Because we could not obtain an analytical expression of $\text{med}(\tau_n)$ in the exponential model, we use the asymptotic approximation given in Equation \eqref{med_asymp}, which is in fact extremely accurate, even for small sample sizes. Thus, the bounds $a_n$ and $b_n$ used for the classification are
\begin{equation}\label{anbn} a_n = \frac 1 {\log(n)} + \frac{\log \log(2)}{\log^2(n)} \ \text{ and } \ b_n = \frac{n+1}{2n}. \end{equation}

Although it is quite straight-forward to implement, the proposed procedure performs well in terms of classification compared to other estimation methods such as Zhang and Stephens (ZSE) \cite{zhang2009new}, Song and Song (SSE) \cite{song2012quantile} or maximum likelihood (MLE). A comparative study is shown in Table \ref{tab:compare}.

\begin{table}[H]
{\small
\begin{center}
\resizebox{\textwidth}{!}{
\begin{tabular}{l||ccc|ccc|ccc|ccc}
& \multicolumn{3}{c|}{Model A} & \multicolumn{6}{c|}{Model B} & \multicolumn{3}{c}{Model C} \\
& \multicolumn{3}{c|}{$k=-0.1$} & \multicolumn{3}{c|}{$k=0.1$}  & \multicolumn{3}{c|}{$k=0.9$} & \multicolumn{3}{c}{$k=1.1$} \\
& \textbf{A} & B & C & A & \textbf{B} & C & A & \textbf{B} & C & A & B & \textbf{C} \\
\hline
ZSE & & & & & & & & & & & \\
$n=15$ & $\mathbf{65.8}$ & $34.1$ & $0.1$ & $44.1$ & $\mathbf{55.5}$ & $0.4 $ & $0.9$ & $\mathbf{76.5}$ & $22.6$ & $0.4$ & $57.7$ & $\mathbf{41.9}$ \\
$n=100$ & $\mathbf{84.2}$ & $15.8$ & $0.0$ & $19.9$ & $\mathbf{80.1}$ & $0.0 $ & $0.0$ & $\mathbf{90.7}$ & $9.3$ & $0.0$ & $31.0$ & $\mathbf{69.0}$ \\
\hline
SSE & & & & & & & & & & & \\
$n=15$ & $\mathbf{85.3}$ & $14.5$ & $0.3$ & $73.8$ & $\mathbf{25.6}$ & $0.7 $ & $20.7$ & $\mathbf{50.0}$ & $29.4$ & $14.7$ & $40.7$ & $\mathbf{44.6}$ \\
$n=100$ & $\mathbf{70.1}$ & $30.0$ & $0.0$ & $29.0$ & $\mathbf{70.9}$ & $0.0 $ & $0.0$ & $\mathbf{77.1}$ & $22.9$ & $0.0$ & $27.3$ & $\mathbf{72.7}$ \\
\hline
MLE & & & & & & & & & & & \\
$n=15$ & $\mathbf{41.4}$ & $50.2$ & $8.4$ & $18.8$ & $\mathbf{64.6}$ & $16.6 $ & $0.1$ & $\mathbf{11.2}$ & $88.8$ & $0.0$ & $4.1$ & $\mathbf{95.9}$ \\
$n=100$ & $\mathbf{75.0}$ & $25.0$ & $0.0$ & $9.7$ & $\mathbf{90.3}$ & $0.0 $ & $0.0$ & $\mathbf{49.8}$ & $50.2$ & $0.0$ & $4.4$ & $\mathbf{95.6}$ \\
\hline
\hline
 & & & & & & & & & & & \\
$n=15$ & $\mathbf{59.7}$ & $39.9$ & $0.4$ & $38.6$ & $\mathbf{60.3}$ & $1.1 $ & $0.2$ & $\mathbf{59.6}$ & $40.2$ & $0.0$ & $40.7$ & $\mathbf{59.3}$ \\
$n=100$ & $\mathbf{74.4}$ & $25.6$ & $0.0$ & $20.4$ & $\mathbf{79.6}$ & $0.0 $ & $0.0$ & $\mathbf{79.8}$ & $20.2$ & $0.0$ & $22.3$ & $\mathbf{77.7}$ \\

\end{tabular}}
\label{tab:compare}
\caption{\footnotesize{Percentages of the model classifications obtained by Zhang and Stephens Estimation (ZSE), Song and Song Estimation (SSE) and Maximum Likelihood Estimation (MLE) compared to the proposed methodology (bottom lines). In bold are the percentages of correct classification.}}
\end{center}}
\end{table}
\vspace*{-0.5cm}

\begin{figure}[H]
    \centering
      \hspace*{-0.1cm}  \includegraphics[width=0.49\textwidth]{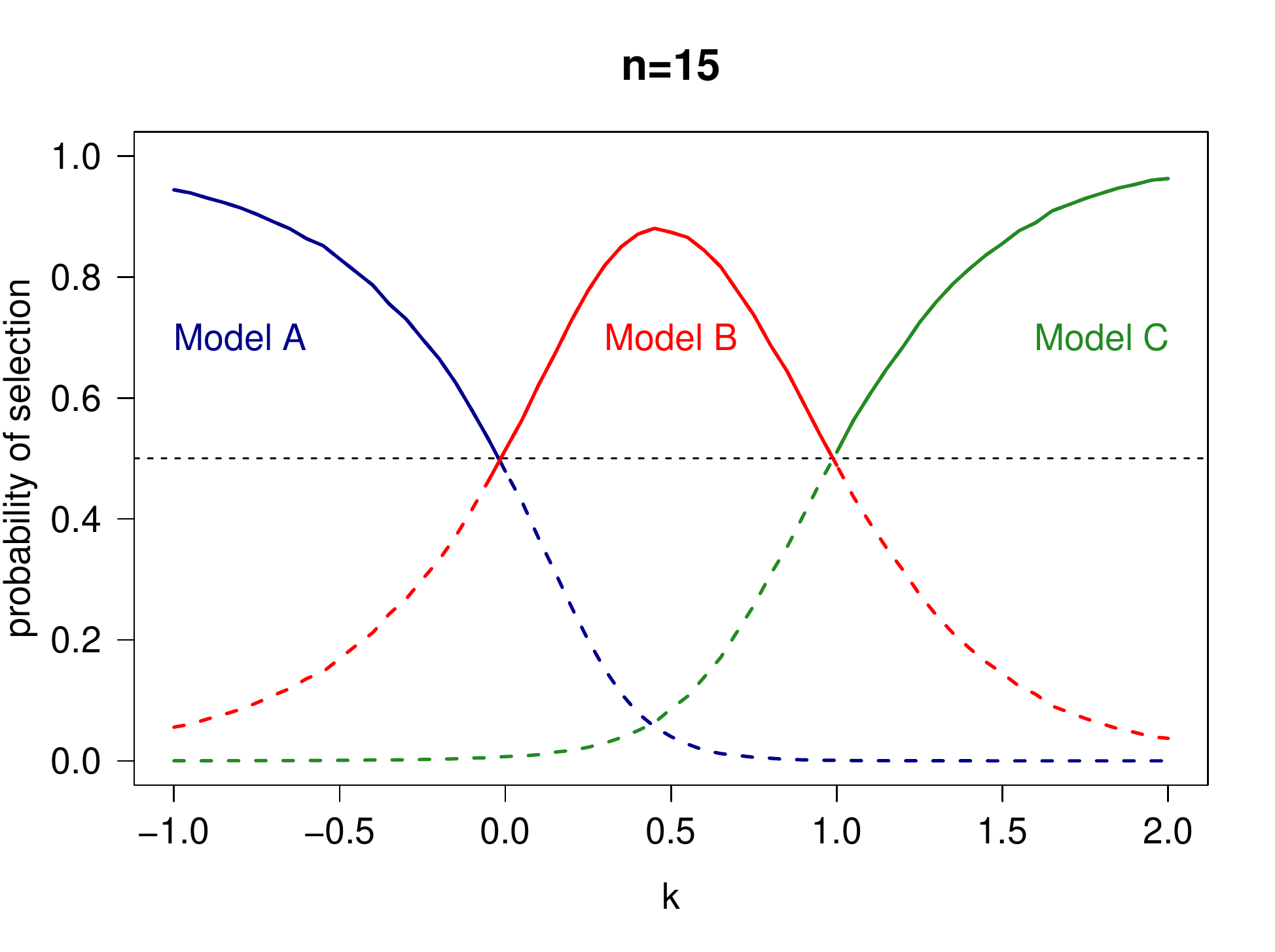}  \hspace*{-0.1cm}  \includegraphics[width=0.49\textwidth]{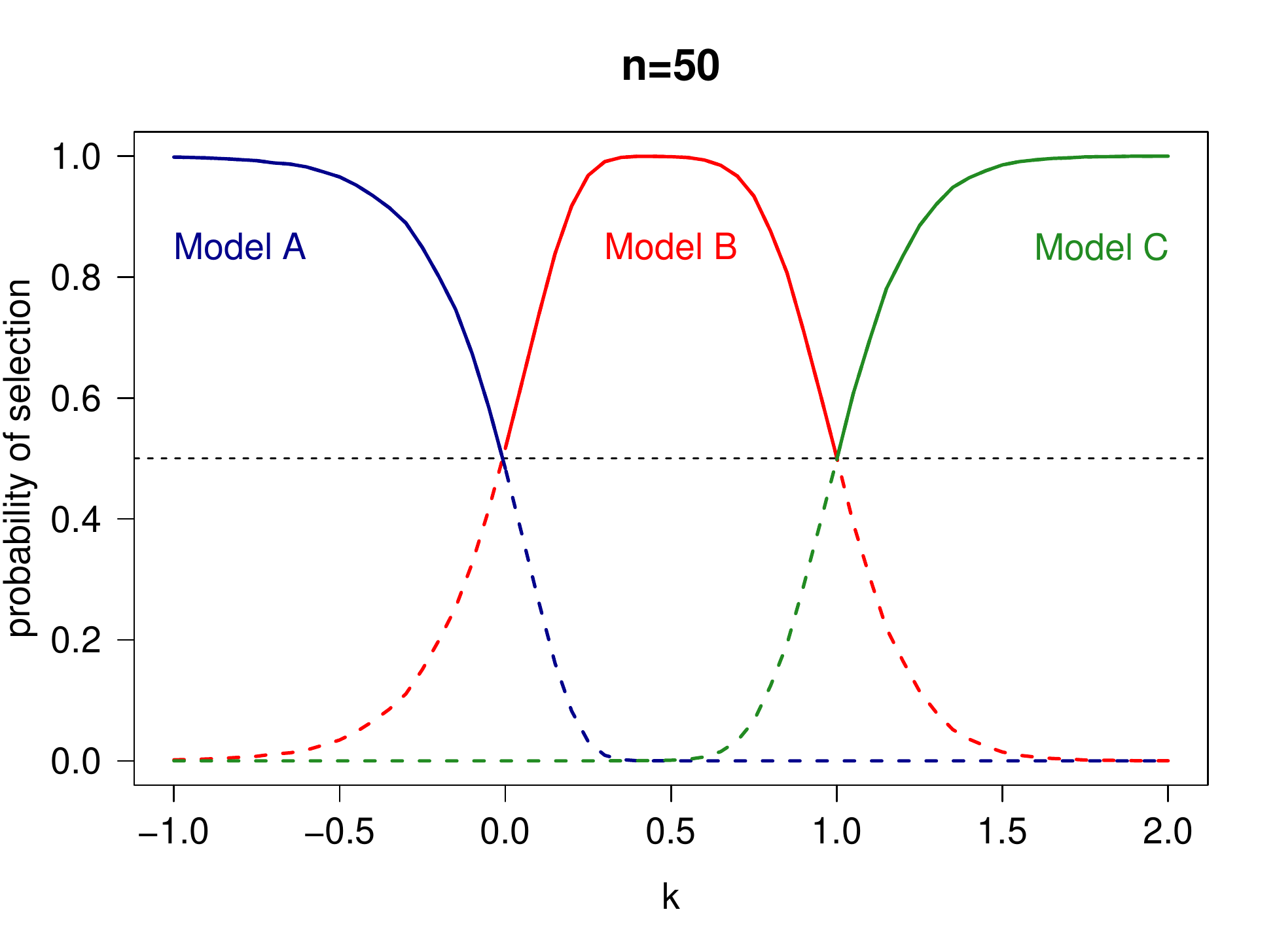}
    \caption{\footnotesize{Monte-Carlo estimated probabilities of selecting the models A,B and C in function of the shape parameter $k$ of the GPD, for sample sizes of $n=15$ (left) and $n=50$ (right). The straight line gives the probability of selecting the correct model. }}
    \label{fig:acc_xi}
\end{figure}

Our classification method manages to not favor any model in particular, as the probability of correct selection depends essentially on the distance between the scale parameter $k$ and the closest transitional phase. This is well illustrated in Figure \ref{fig:acc_xi} where the probability of accurate selection appears to be nearly symmetric locally around $0$ and $1$. In comparison, classification by MLE has a clear tendency to overestimate the scale parameter while SSE tends to underestimate it. For a small sample size, these biases are particularly clear for $k=0.9$ where model C is selected $88.8 \%$ of the time and for $k=0.1$ where model A is selected $73.8 \%$ of the time by SSE. In fact, for the small sample case $n=15$, our method is almost objectively best compared to the three other methods. For the examples considered, the situation $n=100, k \approx 0$ is the only scenario where our method is outperformed, ZSE being clearly best. This strongly suggests that the information on the tail of the distribution is not contained entirely in the mean/max statistic in this case. Nevertheless, the results of the method are overall quite satisfying given the small computational cost. \vspace*{-0.5cm}

\section{Applications}\label{sec:appli}

\subsection{Fit a Generalized Pareto distribution}

Due to the missclassification issue in Generalized Pareto models, inference methods tend to specialize in a specific region of the set of parameters, with none being uniformly best \cite{del2015likelihood}. To solve the problem, we propose a simple recipe that involves a two-step procedure. We assume that the data $X=(X_1,\dots ,X_n)$ are independent and identically distributed from a GPD and that no prior information is available on the parameters. \\

\textbf{Step 1.} Compute the mean/max statistic $\tau_n = \overline X_n/ X_{(n)}$ and compare it to the values $a_n$ and $b_n$ obtained in Eq. \eqref{anbn}. Select:
    \begin{itemize}
        \item Model A if $\tau_n < a_n$
        \item Model B if $a_n< \tau_n \leq b_n$
        \item Model C if $\tau_n \geq b_n$
    \end{itemize}

\textbf{Step 2.} Consider the estimation method depending on the selected model and sample size by this recommendation based on the analysis in \cite{del2015likelihood}. \\

    \begin{table}[htbp]
        \centering
            \begin{tabular}{|c||c|c|c|}
            \hline
            sample size & Model A & Model B & Model C \\
            \hline
            $n \leq 30$ & SSE & ZSE & MLE \\
            \hline
            $n > 30$ & ZSE & ZSE & MLE \\
            \hline

            \end{tabular}
    \end{table}

\noindent In all situations, the method must be applied with the corresponding restriction on the set of parameters. Furthermore, when the sample is classified as Model B or C, a previous estimation of the endpoint can be recommended. Non-parametric approach shows satisfactory results, see \cite{FragaAlves2016}.\\

\begin{figure}[htbp]
        \centering
            \hspace*{0.05cm}        \includegraphics[width=0.49\textwidth]{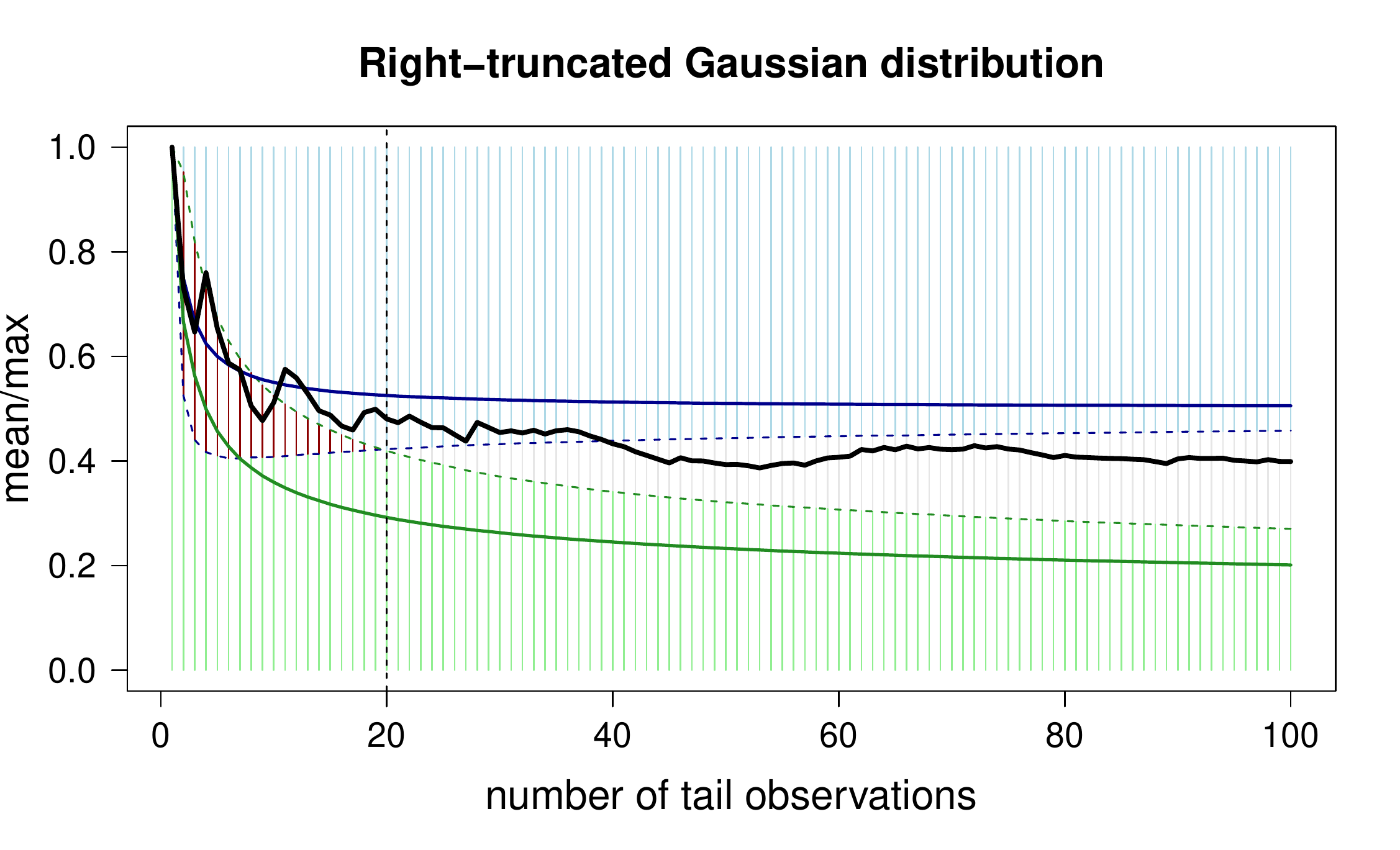}   \hspace*{-0.1cm} \includegraphics[width=0.49\textwidth]{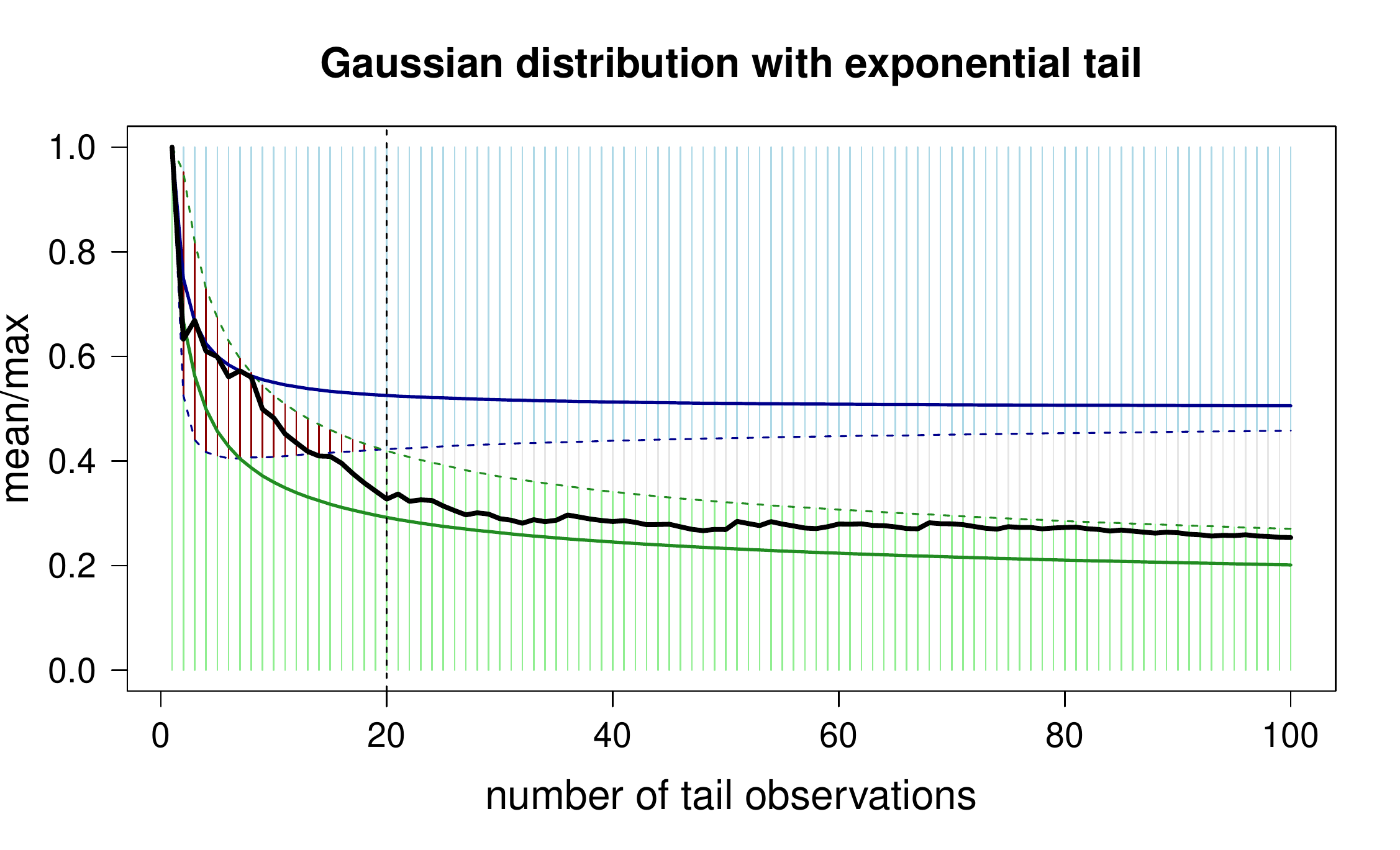}
        \label{fig:ex_gauss}
        \caption{\footnotesize{Examples of trajectories of the mean/max statistic $\tau_n$ from a sample of size $100$ drawn from a standard Gaussian variable with modified tail. On the left, the distribution is truncated at the value $1.5$ (left), on the right, it is extended to an exponential tail starting from $1.5$. The bounds $a_n,b_n$ from  Eq. \eqref{anbn} are represented in straight lines along with their respective $95\%$ confidence regions. The red-striped area corresponds to potential values of $\tau_n$ with more than $5 \%$ uncertainty to distinguish between uniform and exponential tails.}}
    \end{figure}

This recipe is meant for data independently drawn from a Generalized Pareto distribution. For general problems related to tail detection and calibration, this assumes implicitly that the observations have been already extracted from the tail values of a larger dataset and suitably shifted so as to fit the GPD corresponding to the tail distribution. Nevertheless, the mean/max statistic can be used directly for tail detection. Suppose one is interested in the tail of the distribution of a sample $Y_1,...,Y_N$, the idea is to calculate the mean/max statistic over the shifted tail observations $X_i :=Y_{(N-n+i)}- Y_{(N-n)}, i=1,\dots ,n$ for various values of $n$. This way, one obtains a trajectory
$$  \tau_n = \frac{\overline X_n}{X_{(n)}} =  \frac{ \sum_{i=1}^n (Y_{(N-n+i)}- Y_{(N-n)})}{n (Y_{(N)} - Y_{(N-n)})} , n=1,2,... $$
that can be compared to the corresponding values of the bounds $a_n,b_n$. Ideally, only the values in the tail of the distribution are to be used for the classification so that $n$ need not be  too large. The optimal threshold for tail detection may appear more or less clearly (see Figure \ref{fig:ex_gauss}).

        %

    An  overview of the mean/max trajectory can be helpful in order to select the correct model for the tail. In the examples displayed  in Figure \ref{fig:ex_gauss}, the  criterion is most conclusive around $n=20$ even though a sample a this size may contain observations that do not obey the tail distribution. As a possible rule of thumb, the model for the tail can be decided based on the last value of $\tau_n$ that leaves the red-striped area, thus aiming for a precision of over $95 \%$.


\subsection{Distribution of extreme seismic events}

The Gutenberg-Richter (GR) law states that the distribution of seismic moment corresponds a to power-law distribution \cite{Knopoff_Kagan77,Corral_Lacidogna} with a density at a value $M$ proportional to $1/M^{1+\beta}$ for a $\beta$ evaluated to be approximately equal to $0.65$.
However, it is suspected that the collected data on seismic moments show a deviation of GR law for large values of $M$, entailing that the power-law model might have to be extended in order to add an exponential decay above a certain threshold \cite{Corral_Lacidogna}. In \cite{Kagan_gji02}, Kagan enumerates the requirements that an extension of the GR law should fulfill. He also argues however that available seismic catalogs do not allow the reliable estimation of the threshold, except in the global case where the faster exponential decay may take place at the highest observed values of $M$, for which the available data are very poor \cite{Zoller_grl}.  \\

The existence of a theoretical maximal magnitude earthquake entails that the true distribution of seismic moments must be right-truncated, thus displaying a uniform tail, see \cite{del2015likelihood}. Nevertheless, the question remains as to know if the largest seismic events on record are not best modeled by an exponential decay. In order to infer on the underlying physical model, we compute the max/mean plot of the seismic moment collected data from the centroid moment tensor (CMT) catalog \cite{Ekstrom2012}. The analysis is restricted to shallow events (depth $ < 70 $ km) from 1976 to 2013, as recommended in \cite{Kagan_gji02}. The results are displayed in Figure \ref{fig:seismic} below.

\begin{figure}[H]
    \centering
        \includegraphics[width=0.7\textwidth]{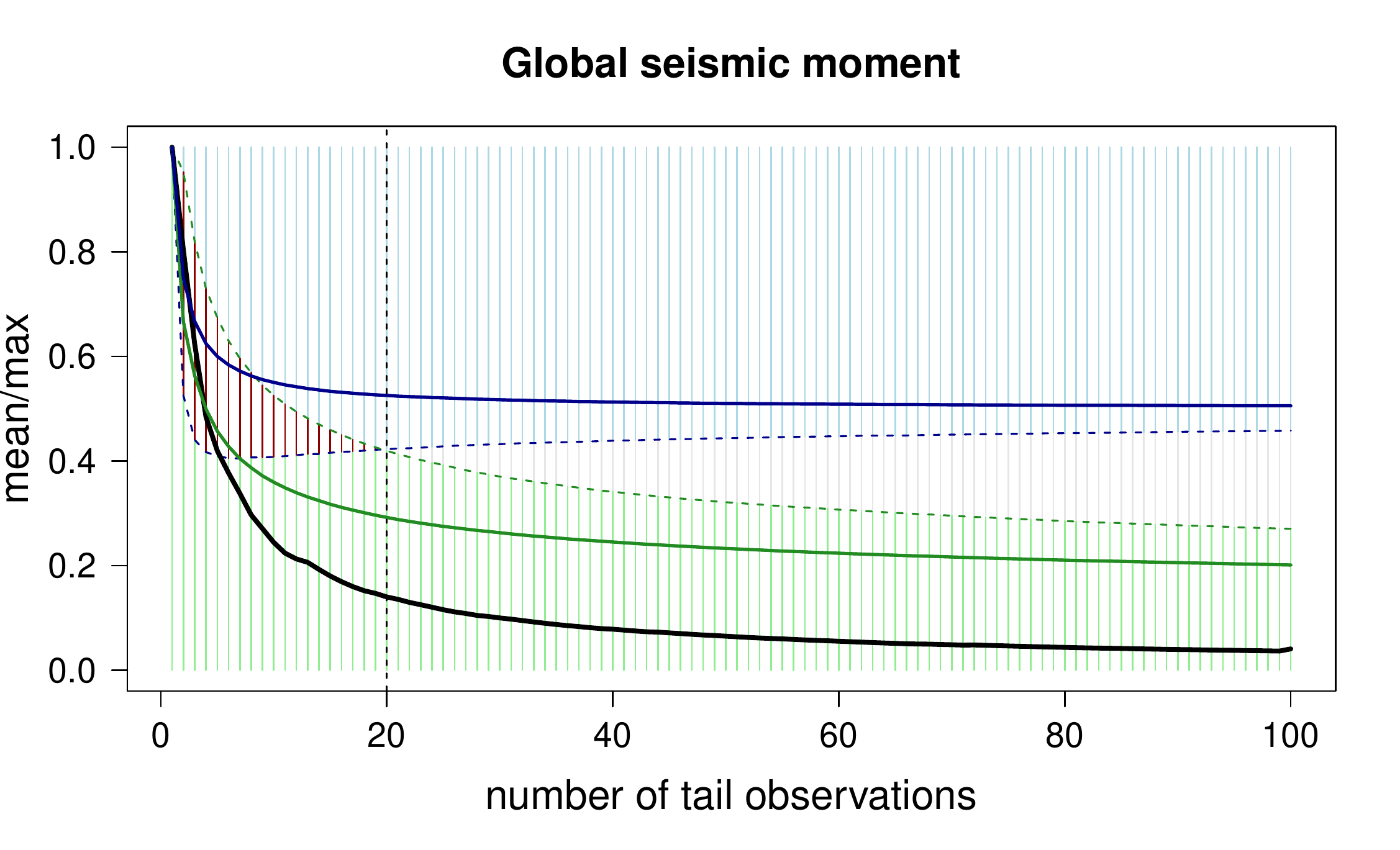}
    \caption{\footnotesize{Trajectory of the mean/max statistic computed over the $n$ largest  (shifted) global seismic events, for $n=1$ to $n=100$.}}
    \label{fig:seismic}
\end{figure}

This  analysis  shows no reason to suspect a change in the tail of the distribution, and confirms that the tail must be classified as Model A with $k < 0$. Consequently, the existence of a cut-off close to the maximal observed seismic moments has to be discarded.

\subsection{Detection of experimental limitations}

\begin{figure}[H]
    \centering
        \includegraphics[width=\textwidth]{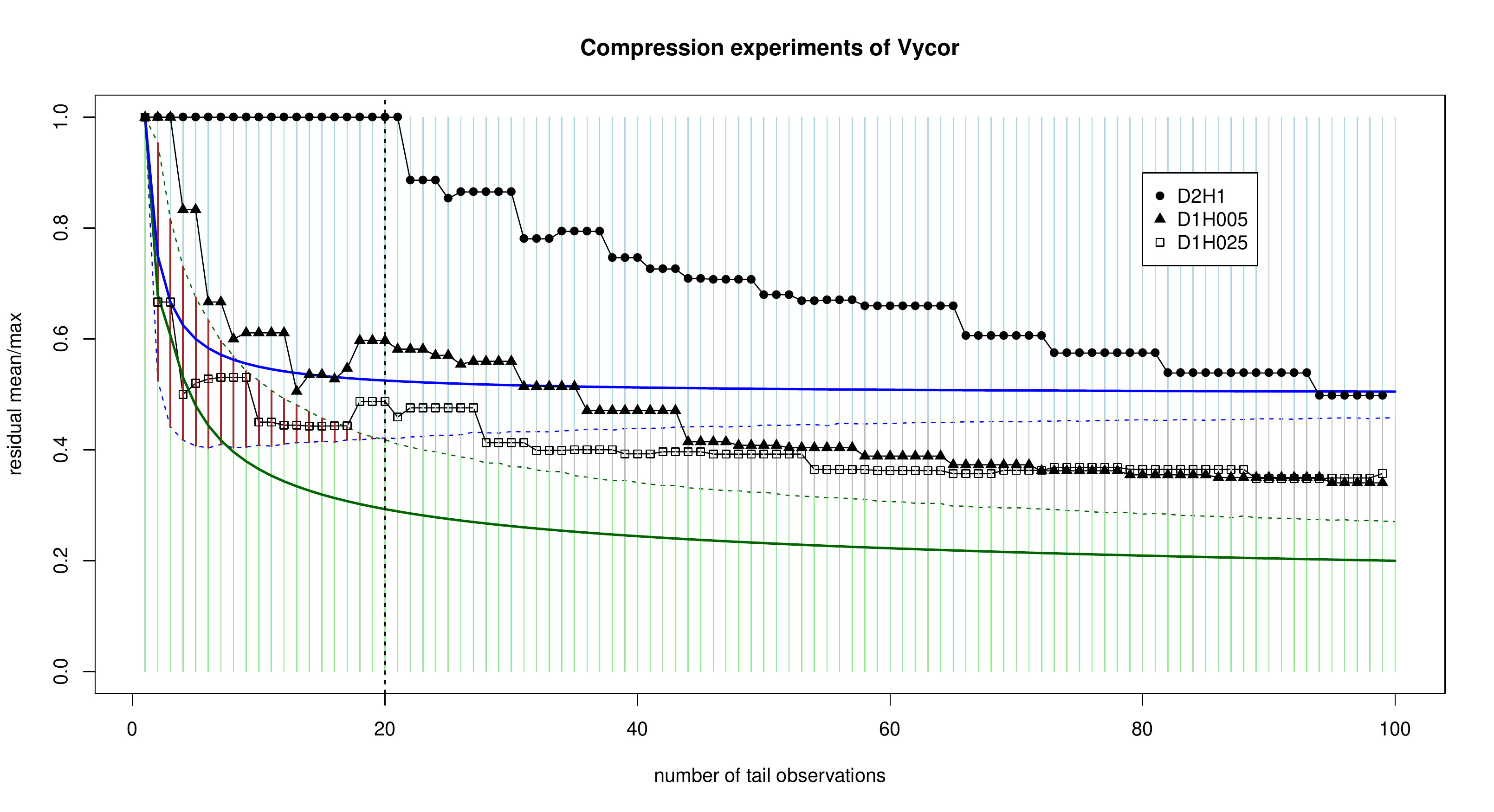}
\caption{\footnotesize{Trajectories of the mean/max statistic for
the 100 largest values obtained from three compression experiments
of nanoporous quartz samples.}}
    \label{fig:ex_vycor}
\end{figure}

We apply this methodology in the analysis of acoustic emission data
in failure under compression experiments of a nanoporous quartz
sample \cite{BaroPRL,NavasCorralVives}. A cylindrical sample is
placed between to plates that approach each other at a constant
velocity. Compression is done with no lateral confinement and the
experiment finishes when the sample has experienced a big failure
and has literally disappeared. In this compressive process, two
transductors are embedded in both plates in order to detect acoustic
emission activity. These signals are pre-amplified in order to
measure properly different magnitudes: amplitude, duration, energy,
etc. The discrete measure in dB of the signal amplitude is defined
as $A=20 \log_{10} (A_{v}/1 \mu V)$, where $A_{v}$ is the highest
voltage value of each acoustic emission signal and $1\mu V$ is a
reference voltage. Remark that the exponential distribution for the
tail, the largest values, is the most natural model from a physical
point of view. Regarding this physical magnitude, one must take into
account that the signal pre-amplification can lead to a saturation
for the largest acoustic emission events. Therefore, this saturation
must be understood as a fictitious cut-off which is inherent to the
experimental set-up. Actually, this experimental limitations are
present in all measuring devices since all of them have a certain
measurement range. The mean/max statistic $\tau_n$ enables to detect
this experimental artifact in most samples, as we see Figure
\ref{fig:ex_vycor} where V.Navas and E.Vives (private communication)
experiments are shown.

According to the experimental intuition, one expects that acoustic emission signals
do not reach saturation if the contact area between
the plates is small. The mean/max methodology provides a
simple and automatic way to identify if the data extracted from any
experimental measurement are saturated or not. Figure
\ref{fig:ex_vycor} shows the trajectories of the mean/max statistic
for the 100 largest values obtained from three compression
experiments of nanoporous quartz samples. The circle trajectories
correspond to the greatest material, 2mm diameter. This mean/max
plot shows saturated experiment, since for the 20 largest emission
signals the exponential decay is rejected. In
fact, the tail of the dataset is classified as a compact support, Model C.
In this case, a descriptive analysis reveals a clear suspicion of
saturation experiment, since the largest value is repeated too many
times. However, for small materials, 1mm diameter, it is difficult to
detect if the saturation occurs. The mean/max test (for the 20
largest) reveals that the triangle trajectories was affected by the
saturation, as in the previous case, but the trajectory of squares
is classified as Model B, therefore a non-saturation can not be
rejected, since Model B includes the exponential distribution for
the tail.

\section*{Acknowledgement}

We are grateful to V. Navas and E. Vives for their feedback.
Research expenses were founded by projects MTM2012-31118 from
Spanish MINECO, 2014SGR-1307 from AGAUR, and the Collaborative
Mathematics Project from La Caixa Foundation.

\appendix

\section{Distribution of the mean/max statistic}

We investigate the distribution of the mean/max statistic $\tau_n$ when the sample is drawn from a GPD. The two transitional phases corresponding to uniform and exponential distributions are given a particular interest as they are directly linked to the level and power of the test discussed in Section \ref{sec:test}.

\subsection{The uniform case}

The threshold $c_\alpha$ corresponding to a test of significance level $\alpha \in (0,1)$ in Theorem \ref{th1} follows from calculating the distribution of $\tau_n$ under the null hypothesis of a uniformly distributed sample $X_1,\dots ,X_n$. Noticing that
$$ \tau_n = \frac 1 n \Big( 1 +  \sum_{i=1}^{n-1} \frac{X_{(i)}}{X_{(n)}} \Big),  $$
it appears that $n(\tau_n- 1)$ is distributed as the sum of $n-1$ independent standard uniform random variables. This distribution is known as the Irwin-Hall distribution of parameter $n-1$, whose density is given by
$$  I_{n-1}(t) =  \frac n {2} \sum_{k=0}^{n-1}  \sgn(t-k) \frac{(-1)^k (t-k)^{n-2}}{k! (n-k)!}, \ t \in [0,n-1],  $$
where $\sgn(y)$ is the sign of $y$, equal to $1$ if $y$ is positive, $-1$ if $y$ is negative and $0$ if $y$ is zero. Geometrically, $I_{n-1}(t)$ represents the volume of the intersection of the $\ell^1$-sphere of radius $t$ in $\mathbb R^{n-1}$ with the hypercube $[0,1]^{n-1}$. For sake of completeness, we derive the distribution of $\tau_n$ under the null hypothesis in the next lemma. The density of a random variable $Z$ will be denoted by $f_Z$.

\begin{lemma}\label{ih} If $X=(X_1,\dots ,X_n)$ is an iid sample from a uniform distribution on $[0,\theta]$, with $\theta >0$, then $\tau_n$ has density
$$ f_{\tau_n}(t) = n I_{n-1} ( n t-1 ) \ , \ t \in \Big[ \frac 1 n,1 \Big].  $$
In particular, $\mathbb E(\tau_n) = \operatorname{med}(\tau_n) = \dfrac{n+1}{2n}$.
\end{lemma}

\noindent \textit{Proof.} First remark that the distribution of $\tau_n$ does not depend on $\theta$ so that we can assume that $\theta=1$ without loss of generality. Let $W_1=X_{(1)}$, $W_i = X_{(i)} - X_{(i-1)}$ for $i=2,\dots ,n$ and $W_{n+1} = 1-X_{(n)}$. We use that $W=(W_1,\dots ,W_{n+1})$ has a standard Dirichlet distribution Dir$(n+1)$, which means that
$$ W \overset{d}{=} \bigg(\frac{Y_1}{\sum_{i=1}^{n+1} Y_i},\dots , \frac{Y_{n+1}}{\sum_{i=1}^{n+1} Y_i} \bigg)  $$
where $Y_1,\dots ,Y_{n+1}$ is an iid sample with exponential distribution. Thus, the vector
$$ W'= \bigg(\frac{X_{(1)}}{X_{(n)}},\dots , \frac{X_{(n-1)}- X_{(n-2)}}{X_{(n)}} \bigg) \overset{d}{=} \left(\frac{Y_1}{\sum_{i=1}^{n} Y_i},\dots , \frac{Y_{n}}{\sum_{i=1}^{n} Y_i} \right) $$
has Dir$(n)$ distribution and $Z_i = X_{(i)}/X_{(n)} \overset{d}{=} \sum_{j=1}^i Y_j/\sum_{j=1}^n Y_j, i =1,\dots ,n-1 $ has the same distribution as an ordered iid sample of $n-1$ uniform random variables on $[0,1]$. We deduce that $\sum_{i=1}^{n-1} Z_i = n(\tau_n-1)$ has Hirwin-Hall distribution with parameter $n-1$ and the result follows. \qed

\subsection{The exponential case}

Calculating the power of the test requires to know the distribution of $\tau_n$ under $\mathcal H_1$. Here again, this distribution can be computed explicitly.

\begin{lemma}\label{exp} If $X=(X_1,\dots ,X_n)$ is an iid sample from an exponential distribution of parameter $\lambda >0$, then $\tau_n$ has density
$$ f_{\tau_n}(t)  = \frac{n!}{n^{n-1}}  \ \frac{ I_{n-1}(n t-1)}{t^n} , \ t \in \Big[ \frac 1 n,1 \Big].   $$
\end{lemma}

\noindent \textit{Proof.} The distribution of $\tau_n$ does not depend on $\lambda$ so that we can assume that $\lambda=1$ without loss of generality. We know that the ordered sample has density on $\mathbb R^n$ given by
$$ (X_{(1)},\dots ,X_{(n)}) \sim  n! \ \exp \big(- \textstyle \sum_{i=1}^n x_i \big) \mathds 1 \{ 0 \leq x_1 \leq \dots  \leq x_n \}.   $$
Let $Y_i = X_{(1)}/X_{(n)}$ for $i=1,\dots ,n-1$. By the change of variable $y_i = x_i/x_n$, we get
$$ (Y_1,\dots ,Y_{n-1},X_{(n)}) \sim n! \ x_n^{n-1} \exp \big(- \textstyle x_n (1+\sum_{i=1}^{n-1} y_i) \big) \mathds 1 \{ 0 \leq y_1 \leq \dots  \leq y_{n-1} \leq 1 ,  x_n > 0 \} .  $$
Integrating the density over $x_n$, we find
$$ Y = (Y_1,\dots ,Y_{n-1}) \sim  \frac{n! (n-1)!}{(1 +\sum_{i=1}^{n-1} y_i)^n} \ \mathds 1 \{ 0 \leq y_1 \leq \dots  \leq y_{n-1} \leq 1 \}   $$
To compute the density of $S_n:= n \tau_n -1= \sum_{i=1}^{n-1} Y_i$, we now need to integrate the joint density over the level sets of the $\ell^1$-norm $D_s = \{ y \in \mathbb R^{n-1}: \sum_{i=1}^{n-1} y_i = s \}$. We obtain
\begin{align} f_{S_n}(s) & = \int_{D_s}  \frac{n! (n-1)!}{(1 +\sum_{i=1}^{n-1} y_i)^n} \ \mathds 1 \{ 0 \leq y_1 \leq \dots  \leq y_{n-1} \leq 1 \} \ dy_1 \dots  dy_{n-1} \nonumber \\
& = \frac{n!}{(1 + s)^n} \int_{D_s}  (n-1)! \ \mathds 1 \{ 0 \leq y_1 \leq \dots  \leq y_{n-1} \leq 1 \} \ dy_1\dots dy_{n-1}. \nonumber
\end{align}
Remark that $(n-1)! \ \mathds 1 \{ 0 \leq y_1 \leq \dots  \leq y_{n-1} \leq 1 \}$ is the density of an ordered sample $(U_{(1)},\dots ,U_{(n-1)})$ of independent uniform random variables on $[0,1]$. Thus,
$$\int_{D_s}  (n-1)! \ \mathds 1 \{ 0 \leq y_1 \leq \dots  \leq y_{n-1} \leq 1 \} \ dy_1\dots dy_{n-1} = I_{n-1}(s)$$
to which we deduce that $S_n=n \tau_n-1$ has density $n! I_{n-1}(s)/(1+s)^n$ over $[0,n-1]$. The result follows by a simple change of variable. \qed

\subsection{Asymptotic distribution of $\tau_n$ in the GPD model}

Finally, we investigate the distribution of $\tau_n$ when the observations are drawn from generalized Pareto distributions with shape parameter $k$. Because $\tau_n$ is scale-free, the following results do not depend on the scale parameter $\sigma$ of the GPD, which can be taken equal to one without loss of generality. The actual distribution of $\tau_n$ being difficult to compute as a function of $k$, we only discuss the asymptotic distribution.

\begin{theorem}\label{gpd_asymp}
Let $X=(X_1,\dots ,X_n)$ be an iid sample drawn from a Generalized Pareto distribution with scale parameter $k$. Then,

\begin{itemize}
    \item If $k > 0$,
    $$ \tau_n \sim \frac k {k+1} + \frac{N_k}{\sqrt n} + o_P \Big( \frac 1 {\sqrt n} \Big),   $$
    where $N_k$ has normal distribution $\mathcal N \Big( 0, \dfrac{k^2}{(1+k)^2 (1+2k)} \Big)$.
    \item If $k=0$ (exponential case),
    $$  \tau_n \sim \frac 1 {\log (n)} - \frac{G}{\log ^2(n)} + o_P \Big( \frac 1 {\log ^2 (n)} \Big),   $$
    where $G$ has standard Gumbel distribution.
    \item If $-1 < k < 0$,
    $$ \tau_n \sim \frac{W_k}{n^{-k}} + o_P(n^{k})$$
    where $W_k$ has Weibull distribution with shape parameter $-1/k$ and scale parameter $-k/(k+1)$.
\end{itemize}
\end{theorem}

\noindent \textit{Proof.} We tackle each case separately. For $k >0$, we have by the central limit theorem,
$$ \sqrt n \Big(\overline X_n - \frac 1 {1+k} \Big) \xrightarrow[n \to \infty]{d} \mathcal N \Big(0, \dfrac{1}{(1+k)^2 (1+2k)} \Big) $$
while $X_{(n)}$ converges a.s.~towards $1/k$, the upper bound of the support of the GPD. The result follows easily by Slutsky's lemma. In the exponential case $k=0$, $\overline X_n$ converges a.s.~towards $1$ while for $t>-\log(n)$,
$$ \mathbb P \Big(X_{(n)} -  \log (n) \leq t  \Big) = \Big( 1 -e^{- (\log (n) + t) } \Big)^n = \Big( 1 - \frac 1 {n^{ 1 + \frac t {\log(n)}}} \Big)^n \xrightarrow[n \to \infty]{} e^{-e^{-t}}.  $$
Thus $X_{(n)} - \log(n)$ converges towards a standard Gumbel distribution. The result follows in view of
$$ \frac 1 {X_{(n)}} = \frac 1 {\log(n) - (X_{(n)} - \log(n))} =\frac 1 {\log(n)} + \frac {X_{(n)} - \log(n)}{\log^2(n)} + o_P \Big( \frac 1 {\log ^2 (n)} \Big).  $$
For the final case $-1 < k < 0$, we have by the strong law of large numbers
$$ \overline X_n = \frac 1 {1+k} + o_P(1) $$
For $t>n^k$, write
$$ \mathbb P \Big( n^{k} (1- k X_{(n)}) \leq t \Big) =\mathbb P \Big( X_{(n)} \leq \frac 1 k \big( 1 - \frac t {n^k} \big) \Big) = \Big( 1 - \frac{t^{1/k}}{n} \Big)^n \xrightarrow[n \to \infty]{} e^{-t^{1/k}}. $$
We deduce that $n^k(1-kX_{(n)}) \sim - k n^k X_{(n)}$ converges towards a Fr\'echet distribution with shape parameter $-1/k$ as $n \to \infty$. Thus, $ -n^{-k}/kX_{(n)}$ converges to the inverse of a Fr\'echet variable whose distribution is Weibull. The result follows by Slutsky's lemma. \qed

\bibliographystyle{plain}
\bibliography{refs}

\end{document}